\documentstyle[12pt]{amsart}
\newtheorem{thm}{Theorem}
\newtheorem{lemma}{Lemma}
\theoremstyle{definition}
\theoremstyle{remark}
\newtheorem{rem}{Remark}
\newtheorem{notation}{Notation}

\newcommand{\M}{{\cal M}}
\newcommand{\Mbar}{{\overline{\cal M}}}
\newcommand{\WP}{{\omega_{WP}}}

\begin{document}
\title[Weil-Petersson volumes]
{Weil-Petersson volumes of moduli spaces\\ of curves
and the genus expansion\\ in two dimensional gravity}
\author{Peter Zograf}
\address[]{Steklov Mathematical Institute, St.Petersburg, 191011 Russia}

\address[]{{\em Current address:} Max-Planck-Institut f\"ur Mathematik, 
53225 Bonn, Germany}
\email[]{zograf@@mpim-bonn.mpg.de}

\date{October 29, 1998}
\thanks{This work was done during the author's stay at the 
Max-Planck-Institut f\"ur Mathematik (Bonn, Germany). The author is grateful 
to MPIM for support and excellent working conditions.}
%\subjclass{Primary ; Secondary }
%\keywords{}
%\dedicatory{}

\begin{abstract}
A formula for the generating function of the
Weil-Petersson volumes of moduli spaces of pointed curves 
that is identical to the genus expansion of the free
energy in two dimensional gravity is obtained. The contribution
of arbitrary genus is expressed in terms of the Bessel function $J_0$.
\end{abstract}

\maketitle

The aim of this note is to prove a rather explicit formula for the generating
function of the Weil-Petersson volumes of the moduli spaces
$\M_{g,n}$ of smooth $n$-pointed curves of genus $g$ that was conjectured
in \cite{Z2}. E.~Getzler noticed that the conecture in \cite{Z2} is identical to  
the formulas in \cite{IZ}, Sect.~5, and \cite{E}, Sect.~3 and Appendix~D, for the genus 
expansion of the free energy in Witten's two dimensional gravity  \cite{G}. 
C.~Faber numerically checked this conjecture for the moduli spaces $\M_{3,n}$ 
with $n\leq 12$ \cite{F}. All that evidence convinced the author that
the conjecture must actually be true not only for the Weil-Petersson volumes
(see Theorem 1 below), but also for their higher analogues in the sense
of \cite{KMZ} (i.e., intesection numbers of Mumford's tautological classes
$\kappa_1,...,\kappa_{3g-3+n}$ on $\Mbar_{g,n}$).

We start by introducing the necessary notation. It is well known that  the K\"ahler form 
$\WP$ of the Weil-Petersson metric on $\M_{g,n}$ extends as a closed
current to the Deligne-Mumford compactification $\Mbar_{g,n}$ 
and represents a real cohomology class $[\WP]\in H^2(\Mbar_{g,n},\bold R)$
which coincides (up to a factor of $\pi^2$) with Mumford's
first tautological class $\kappa_1$ on $\Mbar_{g,n}$ \cite{Wo1, Wo2}. 
It means that the Weil-Petersson volume of the moduli space
$\M_{g,n}$ is finite and is given by
$$Vol_{WP}(\M_{g,n})=\pi^{2(3g-3+n)}\frac{\langle\kappa_1^{3g-3+n}\rangle }
{n!(3g-3+n)!}.$$
(For the definition of intersection numbers
$\langle\kappa_1^{d_1}\kappa_2^{d_2}...\rangle $ we refer to \cite{W}.) 

\begin{notation} Fix $g\geq 0$. The generating function for the
Weil-Petersson volumes of moduli spaces $\M_{g,n}$ is the series
$$\phi_g(x)=\sum_{n=0}^\infty\frac{V_{g,n}}{n!\,(3g-3+n)!}\,x^n,$$
where $V_{g,n}=\langle\kappa_1^{3g-3+n}\rangle \;$ (we assume that
$V_{0,0}=V_{0,1}=V_{0,2}=V_{1,0}=0$ and $V_{0,3}=1$).
\end{notation}

The function $\phi_0(x)$ is completely characterized by the property 
that $y(x)=\phi''_0(x)$ is inverse to the Bessel function 
$x(y)=-\sqrt{y}J'_0(2\sqrt{y})$ \cite{KMZ}. Now, following \cite{IZ} 
(cf. also \cite{E}, Eq.~(2.28)), we introduce a sequence of functions 
$f_1,f_2,...$, where
\begin{equation}
f_1(x)=1-\frac{1}{y'(x)},\quad f_2(x)=\frac{y''(x)}{y'(x)^3},\quad 
f_i(x)=\frac{f'_{i-1}(x)}{y'(x)},\quad i=3,4,... .\label{f}
\end{equation}

\begin{thm} For the generating function $\phi_g$ with $g\geq 2$ we have
\begin{equation}
\phi_g= \sum_{|l|=3g-3}\langle\tau_2^{l_2}\tau_3^{l_3}...\tau_{3g-2}^{l_{3g-2}}\rangle \;
{y'}^{2(g-1)+\|l\|}\; \prod_{i=2}^{3g-2}\frac{f_i^{l_i}}{l_i!}, \label{gen}
\end{equation}
where $l=(l_2,l_3,...,l_{3g-2})$ is a multi-index with $l_i\geq 0$,
$$|l|=\sum_{i=2}^{3g-2}(i-1)\,l_i,\quad\quad \|l\|=\sum_{i=2}^{3g-2}l_i,$$
and the intersection numbers 
$\langle\tau_2^{l_2}\tau_3^{l_3}...\tau_{3g-2}^{l_{3g-2}}\rangle $
are defined in \cite{W}. 
\end{thm}

\begin{rem}
Formula (\ref{gen}) is identical to formula (5.27) in \cite{IZ} that expresses
the genus $g$ contribution $F_g(t_0,t_1,...)$ to the free energy in terms
of the function $u_0=\frac{\partial^2}{\partial t_0^2}\,F_0$. In our setting,
$y(x)$ plays the role of $u_0(t_0,t_1,...)$. Moreover, as it was noticed by
B.~Dubrovin \cite{D},
$$y(x)=u_0(t_0,t_1,...) \left|_{\,t_0=x,\,t_1=0,\,t_k=\frac{(-1)^k}{(k-1)!}\;\;(k=2,3,...)}\right. $$
and, therefore,
$$\phi_g(x)=F_g(t_0,t_1,...) \left|_{\,t_0=x,\,t_1=0,\,t_k=\frac{(-1)^k}{(k-1)!}\;\;(k=2,3,...)}\right. .$$
\end{rem}

The following observation, which we owe to B. Dubrovin (see the remark above), 
will be useful in proving the theorem:

\begin{lemma} The functions $f_i$ defined by relation (\ref{f}) satisfy the 
functional equations 
\begin{equation}
f_i(x)=\sum_{k=0}^\infty \frac{(-1)^{i+k}}{(i+k-1)!}\;\frac{y(x)^k}{k!},
\label{sum}\end{equation}
where $y(x)=\phi''_0(x)$. In particular, for $i\geq 2$
$$f_i(0)=\frac{(-1)^i}{(i-1)!}.$$
\end{lemma}

{\em Proof of Lemma.} According to \cite{KMZ} we have
$$ x=-\sqrt{y(x)}\,J'_0\left(2\sqrt{y(x)}\right)=
\sum_{k=1}^\infty\frac{(-1)^{k-1}}{(k-1)!}\;\frac{y(x)^k}{k!},$$
or, equivalently,
$$y(x)-x=\sum_{k=2}^\infty\frac{(-1)^{k}}{(k-1)!}\;\frac{y(x)^k}{k!}.$$
Applying the differential operator $\left(\frac{1}{y'(x)}\frac{d}{dx}\right)^n$ 
to both sides of the last equality we get the statement of the lemma.

{\em Proof of Theorem.} E. Witten observed in \cite{W} that 
$\kappa$-intersection numbers on moduli spaces of stable curves can
be expressed in terms of $\tau$-intersection numbers and vise versa.
An explicit formula was derived in \cite{KMZ}, Corollary 2.3. In our 
notation, adapted for the Weil-Petersson volume of the moduli space 
$\M_{g,n}$, it reads 
\begin{equation}
\frac{\langle\kappa_1^{3g-3+n}\rangle }{(3g-3+n)!}=\sum_{|l|=3g-3+n}
\langle\tau_0^n\tau_2^{l_2}\tau_3^{l_3}...\tau_{3g-2+n}^{l_{3g-2+n}}\rangle 
\frac{(-1)^{g-1+n+\|l\|}}{\prod_{i=2}^{3g-2+n}l_i!\,((i-1)!)^{l_i}},
\label{vol}\end{equation}
where, as above, $l=(l_2,l_3,...,l_{3g-2+n}),\;l_i\geq 0,\;
\|l\|=l_2+l_3+...+l_{3g-2+n}$ and $|l|=l_2+2l_3+...+(3g-3+n)l_{3g-2+n}$.

Now denote by $\Phi_g$ the right hand side of equation (\ref{gen}). In view
of the above lemma and formula (\ref{vol}),  it suffices to show that for
its $n$-th derivative we have
\begin{equation}
\Phi_g^{(n)}= \sum_{|l|=3g-3+n}\langle\tau_0^n\tau_2^{l_2}\tau_3^{l_3}...
\tau_{3g-2+n}^{l_{3g-2+n}}\rangle \;
{y'}^{2(g-1)+n+\|l\|}\; \prod_{i=2}^{3g-2+n}\frac{f_i^{l_i}}{l_i!}.
\label{der}\end{equation}
We will do that by induction in $n$. The base case $n=0$ being obvious,
let us derive equation (\ref{der}) for $\Phi_g^{(n)}$ assuming that it holds for 
$\Phi_g^{(n-1)}$. Because of the definition of the functions $f_i$, all we have to 
do is to check that
\begin{eqnarray}
\langle\tau_0^n\prod_{i=2}^{3g-2+n}\tau_i^{l_i}\rangle &=&l_2\left(2(g-1)+(n-1)+(\|l\|-1)\right)
\langle\tau_0^{n-1}\prod_{i=2}^{3g-2+n}\tau_i^{l_i-\delta_{i,2}}\rangle \nonumber\\
&+&\sum_{j=3}^{3g-2+n}
l_j\langle\tau_0^{n-1}\prod_{i=2}^{3g-2+n}\tau_i^{l_i-\delta_{i,j}+\delta_{i,j-1}}\rangle .
\label{ind}\end{eqnarray}
Actually, formula (\ref{ind}) is an immediate consequence of the 
puncture and dilaton equations (see \cite{W}, Eqs. (2.41) and (2.45)),
which in our case are, respectively,
$$\langle\tau_0^n\prod_{i=2}^{3g-2+n}\tau_i^{l_i}\rangle \;
=\;\sum_{j=2}^{3g-2+n}l_j\langle\tau_0^{n-1}\tau_1^{\delta_{j,2}}
\prod_{i=2}^{3g-2+n}\tau_i^{l_i-\delta_{i,j}+\delta_{i,j-1}}\rangle $$
and
$$\langle\tau_0^{n-1}\tau_1\prod_{i=2}^{3g-2+n}\tau_i^{l_i-\delta_{i,2}}\rangle \;=\;\left(2(g-1)
+(n-1)+(\|l\|-1)\right)\langle\tau_0^{n-1}\prod_{i=2}^{3g-2+n}\tau_i^{l_i-\delta_{i,2}}\rangle .$$
That completes the proof of the theorem.

\begin{rem}
For $g=2$ the above theorem was proven in \cite{Z2} by a different method
similar to that of \cite{Z1}, which is based on straightforward computations of 
intersections on $\Mbar_{2,n}$ and does not use Witten's theory of two
dimensional gravity.
\end{rem} 

\begin{rem}
Theorem 1 implies that Weil-Petersson volumes behave
asymptotically for large $n$ as
\begin{equation}\label{asympt}
\frac{V_{g,n}}{n!\,(3g-3+n)!}\sim C^{\,n}\,n^{-1+\frac{5}{2}(g-1)},
\end{equation}
with a constant $C$ independent of $g$, as suggested by 
C. Itzykson \cite{I}. For $g=0$ this formula was proven in
\cite{KMZ}, and the constant $C$ was expressed there in terms of the first
zero of the Bessel function $J_0$. 
\end{rem}

\begin{rem}
In fact, the results of \cite{KMZ} combined with the present approach
allow to prove an analogue of Theorem 1 for the generating function
of higher Weil-Petersson volumes, or intersection numbers
$\langle\kappa_1^{d_1}\kappa_2^{d_2}...\rangle$. However, since the
combinatorics becomes more involved and the notation more cumbersome,
the detailed exposition will appear elsewhere.
\end{rem}

\subsection*{Acknowledgment}
It is a pleasure to thank E.~Getzler for bringing the references \cite{IZ}, 
\cite{E} to our attention, C.~Faber for informing us about his numerical
computations, B.~Dubrovin for many useful discussions, and Yu.~Manin
for his encouragement and interest in this work.

\end{document}